\documentclass[fontsize=11pt,a4paper,DIV12]{scrartcl}

\usepackage[utf8]{inputenc} 
\usepackage[boldsans]{concmath} %
\usepackage{euler} %
\usepackage{rsfso}

\usepackage{url}

\usepackage{graphicx}

\def\NI{\noindent}
\def\ni{\noindent}
\def\sk{\smallskip}

\def\Proof{\par{\NI\bf Proof. }}
\def\qed{\hfill\fbox{\hbox{}}\bs}

\def\term#1{{\em #1}\marginpar{\raggedright\textit{\small #1}}}
\usepackage{amssymb}
\usepackage{amsmath}
\usepackage{amsthm}
\usepackage{amsopn}

\usepackage{xfrac}
\usepackage{enumerate}
\usepackage{import}

\newtheoremstyle{meiner} %name
    {4pt}{3pt}           %space above/below
    {\itshape}           %text font
    {}                   %heading indent
    {\sffamily\bfseries} %heading font
    {}                   %punctuation
    { }                  %space after head
    {}                   %head spec ??
\theoremstyle{meiner}

\usepackage[capitalize]{cleveref}

  \newtheorem{thm}{Theorem}
  \crefname{thm}{Theorem}{Theorems}
  \Crefname{thm}{Theorem}{Theorems}
  \newtheorem{prop}{Proposition}
  \newtheorem{lemma}{Lemma}
  \newtheorem{cor}{Corollary}
  
  \crefname{claim}{Claim}{Claims}
  \Crefname{claim}{Claim}{Claims}

  \crefname{obs}{Observation}{Observations}
  \Crefname{obs}{Observation}{Observations}
\theoremstyle{definition}

\newcommand*{\claimproofname}{Proof}

\def\Case#1.{\rih{Case #1.}}
\def\Phase#1.{\NI{\bfPhase #1.}}
\def\Fact#1.{\par\sk{\NI\bf Fact #1.}}
\def\Claim#1.{{\bf Claim #1.}}

\def\Proof{\ni{\slshape Proof.\/}\ }

\def\qed{\hfill\fbox{\hbox{}}\medskip}

\def\ITEMMACRO #1 ??? #2 ???{\par\medskip\noindent%
\hangindent=#2em\setbox0\hbox{#1 \kern5pt}%
\ifdim\wd0<\hangindent\setbox0\hbox to\hangindent{\hss#1\quad}\fi%
\box0\ignorespaces}

\def\Item#1{\ITEMMACRO #1 ??? 2.5 ???}

\def\Bitem{\Item{\hss$\bullet$}}

\def\abstract{\noindent\hfil\vbox\bgroup\hsize=.9\hsize
\small\noindent{\bfseries Abstract.}
}
\def\endabstract{\egroup\hfil}

\graphicspath{{Figures/}}
\input{pdffig.sty}
\usepackage{subfig}
\usepackage{color}
\def\RR{\hbox{\upshape\sffamily I\kern-1ptR}}
\def\NN{\hbox{\upshape\sffamily I\kern-1ptN}}
\def\ZZ{\hbox{\upshape\sffamily Z\kern-4ptZ}}

\newcommand{\PXP}{\mathbf{P}=(X,P)}

\def\A{\mathcal{A}}
\def\C{\mathcal{C}}

\def\mC{C^+}

\def\Min{\mathrm{Min}}

\def\ovl{\overline}
\def\bg{\mathbf{g}}

\usepackage{bibitemsep}

\def\PsFig#1#2#3{}
\def\SC{\scshape}
\begin{document}

\title{4-Connected Triangulations on Few Lines\footnote{%%
  An extended abstract of this work appears in Proceedings of the 27th
  International Symposium on Graph Drawing and Network Visualization
  (GD 19) \texttt{arXiv:1908.04524v1}}}
\author{
\parbox{7.5cm}{\center
{\SC Stefan Felsner%
\footnote{Partially supported by DFG grant FE-340/11-1.}
}\\[3pt]
\normalsize
%%\email{felsner@math.tu-berlin.de}\\
\small
        {Institut f\"ur Mathematik\\
         Technische Universit\"at Berlin}}
}%%END AUTHOR

\date{\vskip-3mm}
\maketitle
%%

% -------------------------------------------------------------
%                   ABSTRACT
% -------------------------------------------------------------
\begin{abstract}
  We show that 4-connected plane triangulations
  can be redrawn such that edges are represented by straight
  segments and the vertices are
  covered by a set of at most $\sqrt{2n}$ lines
  each of them horizontal or vertical.
  The same holds for all subgraphs of
  such triangulations.

  The proof is based on a corresponding result for
  diagrams of planar lattices which makes use of
  orthogonal chain and antichain families. 
\end{abstract}

% -------------------------------------------------------------
%                   SECTION
% -------------------------------------------------------------
\section{Introduction}

Given a planar graph $G$ we denote by $\pi(G)$ the minimum number $\ell$
such that $G$ has a plane straight-line drawing in which the
vertices can be covered by a collection of $\ell$ lines.
Clearly $\pi(G) =1$ if and only if $G$ is a forest of paths.
The set of graphs with $\pi(G)=2$, however, is already
surprisingly rich, it contains trees, outerplanar graphs and
subgraphs of grids, see \cite{BannisterDDEW19,Epp19}.

The parameter $\pi(G)$ has received some attention in recent
years, here is a list of known results:
\Bitem
It is NP-complete to decide whether $\pi(G)=2$
(Biedl et al.~\cite{Berti18}).
\Bitem
For a stacked triangulation $G$, a.k.a.~planar 3-tree or Apollonian network,
let $d_G$ be the stacking depth (e.g.~$K_4$ has stacking depth 1).
On this class lower and upper bounds on $\pi(G)$ are $d_G+1$ and $d_G+2$
respectively, see~Biedl et al.~\cite{Berti18} and for the lower
bound also Eppstein~\cite[Thm. 16.13]{EppBook}.
\Bitem
Eppstein~\cite{Epp19} constructed a planar, cubic, 3-connected,
bipartite graph $G_\ell$ on $O(\ell^3)$ vertices with $\pi(G_\ell) \geq
\ell$.
\smallskip

\ni
Related parameters have been studied by Chaplick et al.~\cite{FLRVW16,FLRVW17}.
\medskip

\begin{figure}[htb]
    \centerline{\input{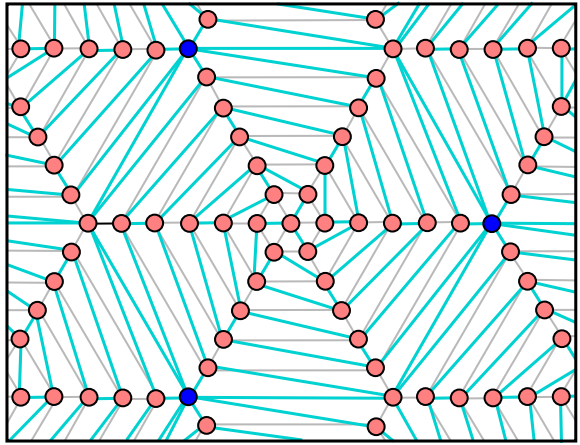}}
    \caption{A section of a 4-connected triangulation
  extending Eppstein's graph $G_\ell$, extra vertices are shown in blue and
  additional edges in cyan.\label{fig:eppstein-extension}}
    \end{figure}

The main result of this paper is the following theorem.
\begin{thm}\label{thm:main}
If $G$ is a 4-connected plane triangulation
on $n$ vertices, then $\pi(G) \leq \sqrt{2n}$.
\end{thm}
The result is not far from optimal since, using a small number of
additional vertices and many additional edges, the graph
$G_\ell$ mentioned above can be transformed into a 
4-connected plane triangulation, i.e.,
in the class we have graphs with $\pi(G) \in \Omega(n^{1/3})$.
Figure~\ref{fig:eppstein-extension} shows an section of such an
extension of $G_\ell$.

%%%%%%%%%%%%%%%%%%%%%%%%%%%%%%%%%%%%%%%%%%%%%%%%%%%%%%%%%%%%%%%%%%%%%% 
% in einem figure environment mit caption
   \calc_figscale{40}
    
%%%%%%%%%%%%%%%%%%%%%%%%%%%%%%%%%%%%%%%%%%%%%%%%%%%%%%%%%%%%%%%%%%%%%%

The proof of the \cref{thm:main} makes use of transversal structures, these
are special colorings of the edges of a 4-connected inner
triangulation of a 4-gon with colors red and blue.

In \cref{ssec:transv-str} we survey transversal structures.  The red
subgraph of a transversal structure can be interpreted as the diagram
of a planar lattice. Background on posets and lattices is given in
\cref{ssec:posets}. Dimension of posets and the connection with
planarity are covered in \cref{subs:dimension}.  In
\cref{ssec:GKtheory} we survey orthogonal partitions of posets.  The
theory implies that every poset on $n$ elements can be covered by
at most $\sqrt{2n}-1$ subsets such that each of the subsets is a chain or an
antichain.

In \cref{sec:plaLat} we prove that the diagram of a
planar lattice on $n$ elements has a straight-line drawing
with vertices placed on a set of $\sqrt{2n}-1$ lines. All the lines
used for the construction are either horizontal or vertical.

Finally in \cref{sec:4ct} we prove the main result: transversal
structures can be drawn on at most $\sqrt{2n}-1$ lines.  In fact, the
red subgraph of the transversal structure has such a drawing by the
result of the previous section.  It is rather easy to add the blue
edges to this drawing. \cref{thm:main} is obtained as a corollary.

%%%%%%%%%%%%%%%%%%%%%%%%%%%%%%%%%%%%%%%%%%%%%%%%%%%%%%
\section{Preliminaries}
%%%%%%%%%%%%%%%%%%%%%%%%%%%%%%%%%%%%%%%%%%%%%%%%%%%%%%
\subsection{Transversal structures}
\label{ssec:transv-str}

Let $G$ be an internally 4-connected inner
triangulation of a 4-gon, in other words
$G$ is a plane graph with quadrangular outer face,
triangular inner faces, and no separating triangle.
Let $s,a,t,b$ be the  outer vertices of $G$ in clockwise order.
A \term{transversal structure}
for $G$ is an orientation and 2-coloring of the
inner edges of $G$ such that

\Item{(1)} All edges incident to $s$, $a$, $t$ and $b$ are 
red outgoing, blue outgoing, red incoming, and blue incoming, respectively.
\Item{(2)} The edges incident to an inner vertex $v$ come in 
clockwise order in four non-empty blocks consisting solely of
red outgoing, blue outgoing, red incoming,  blue incoming edges,
respectively.
\smallskip

\begin{figure}[htb]
    \centerline{\input{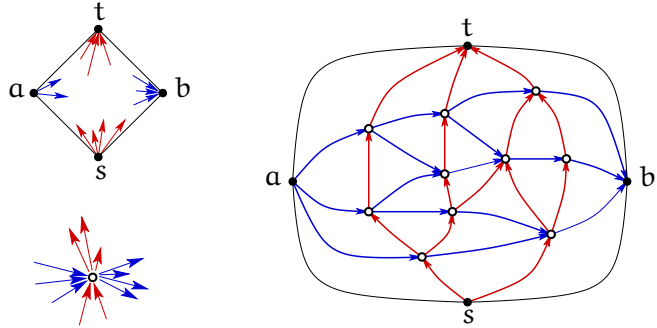}}
    \caption{The two local conditions and an example
  of a transversal structure.\label{fig:trans-struct}}
    \end{figure}

\noindent 
Figure~\ref{fig:trans-struct} illustrates the properties and shows an
example. Transversal structures have been studied
in~\cite{kh-rel4cpga-97}, \cite{f-ccpaa-07}, and \cite{f-tst-09}.
In particular it has been shown that every internally 4-connected inner
triangulation of a 4-gon admits a transversal structure.
Fusy~\cite{f-tst-09} used transversal structures to prove the existence of straight-line
drawings with vertices being placed on integer points $(x,y)$
with $0\leq x \leq W$, $0\leq y \leq H$, and $H+W \leq n-1$. 

%%%%%%%%%%%%%%%%%%%%%%%%%%%%%%%%%%%%%%%%%%%%%%%%%%%%%%%%%%%%%%%%%%%%%% 
% in einem figure environment mit caption
   \calc_figscale{20}
    
%%%%%%%%%%%%%%%%%%%%%%%%%%%%%%%%%%%%%%%%%%%%%%%%%%%%%%%%%%%%%%%%%%%%%%

An orientation of a graph $G$ is said to be acyclic if it has no
directed cycle. Given an acyclic orientation of $G$, a vertex having
no incoming edge is called a \textit{source}, and a vertex having no
outgoing edge is called a \textit{sink}. A \term{bipolar orientation}
is an acyclic orientation with a unique source $s$ and a unique sink
$t$, cf.~\cite{fmr-bor-95}.  Bipolar orientations of
plane graphs are also required to have $s$ and $t$ 
incident to the outer face.  A bipolar orientation of a plane graph
has the property that at each vertex $v$ the outgoing edges form a
contiguous block and the incoming edges form a contiguous
block. Moreover, each face $f$ of $G$ has two special vertices $s_f$
and $t_f$ such that the boundary of $f$ consists of two non-empty
oriented paths from $s_f$ to $t_f$.

Let $G=(V,E)$ be an internally 4-connected inner
triangulation of a 4-gon with outer vertices $s,a,t,b$ 
in clockwise order, and let $E_R$ and $E_B$ respectively be the
red and blue oriented edges of a transversal structure on $G$.
We define $E_R^+ = E_R \cup\{((s,a),(s,b),(a,t),(b,t)\}$
and $E_B^+ = E_B\cup\{((a,s),(a,t),(s,b),(t,b)\}$, i.e.,
we think of the outer edges as having both, a red direction and a
blue direction.
The following has been shown in~\cite{kh-rel4cpga-97} and
in~\cite{f-ccpaa-07}.

\begin{prop}\label{prop:bipolar}
The red graph $G_R=(V,E_R^+)$ and the blue graph $G_B=(V,E_B^+)$ are both
bipolar orientations. $G_R$ has source $s$ and sink $t$, and 
$G_B$ has source $a$ and sink $b$.
\end{prop}

The following two properties are easy consequences of the previous
discussion.

\Item{(R)}
The red and the blue graph are both transitively reduced, i.e.,
if $(v,v')$ is an edge, then there is no directed path
$v,u_1,\ldots,u_k,v'$ with $k\geq 1$.

\Item{(F)}
For every blue edge $e\in E_B$ there is a face $f$ in the red graph 
such that $e$ has one endpoint on each of the two oriented $s_f$ to $t_f$
paths on the boundary of $f$.

%%%%%%%%%%%%%%%%%%%%%%%%%%%%%%%%%%%%%%%%%%%%%%%%%%%%%%
\subsection{Posets}
\label{ssec:posets}

We assume basic familiarity with concepts and terminology for posets,
referring the reader to the monograph~\cite{Trotter-Book} and survey
article~\cite{Trotter-Handbook} for additional background material.
In this paper we consider a poset $P=(X,<)$ as being equipped with a
\textit{strict} partial order.

A \term{cover relation} of $P$ is a pair $(x,y)$ with $x < y$ such
that there is no $z$ with $x < z < y$, we write $x\prec y$ to denote a
cover relation of the two elements.  A \term{diagram} (a.k.a.\ Hasse
diagram) of a poset is an upward drawing of its transitive
reduction. That is, $X$ is represented by a set of points in the plane
and a cover relation $x\prec y$ is represented by a $y$-monotone curve
going upwards from $x$ to~$y$.  In general these curves (edges) may
cross each other but must not touch any vertices other than their
endpoints. A diagram uniquely describes a poset, therefore, we usually
show diagrams in our figures. A poset is said to be \term{planar} if
it has a planar diagram.

It is well known that in discussions of graph planarity, we can
restrict our attention to straight-line drawings.  In fact, using for
example a result of Schnyder~\cite{Schnyder}, if a planar graph has
$n$ vertices, then it admits a planar straight-line drawing with
vertices on an $(n-2)\times(n-2)$ grid.  Discussions of
planarity for posets can also be restricted to straight-line drawings;
however, this may come at some cost in visual clarity. Di Battista et
al.~\cite{diBa-Ta-To} have shown that an exponentially large grid may
be required for upward planar drawings of directed acyclic planar graphs with
straight lines. In the next subsection we will see that for certain
planar posets the situation is more favorable.

%%%%%%%%%%%%%%%%%%%%%%%%%%%%%%%%%%%%%%%%%%%%%%%%%%%%%%
\subsection{Dimension of planar posets}\label{subs:dimension}

Let $P=(X,<)$ be a poset. A \term{realizer} of $P$ is a collection
$L_1,L_2,\dots,L_t$ of linear extensions of $P$ such that
$P=L_1\cap L_2\cap\dots\cap L_t$. The \term{dimension} of $P=(X,<)$,
denoted $\dim(P)$, is the least positive integer $t$ such that $P$ has
a realizer of size~$t$.  Obviously, a poset $P$ has
dimension~$1$ if and only if it is a chain (total order).  Also, there
is an elementary characterization of posets of dimension at most~$2$
that we shall use.

\begin{prop}\label{prop:dim2} 
  A poset $\PXP$ has dimension as most~$2$ if and only if its
  incomparability graph is also a comparability graph.
\end{prop}
%% Proof in posetcover

There are a number of results concerning the dimension of posets with planar
order diagrams.  Recall that an element is called a \textit{zero} of
a poset $P$ when it is the unique minimal element.  Dually, a \textit{one}
is a unique maximal element.  A finite poset which is also a lattice, i.e.,
which has well defined meet and join operations, always
has both a zero and a one.  

The following result  may be considered 
part of the folklore of the subject. 

\begin{thm}\label{thm:2d-lat}
Let $P$ be a finite lattice.  Then $P$ is planar if and
only if it has dimension at most~$2$. 
\end{thm}
%% Proof in posetcover
%% dort fehlt allerdings, dass die Zeichnung Kreuzungsfrei ist

%%%%%%%%%%%%%%%%%%%%%%%%%%%%%%%%%%%%%%%%%%%%%%%%%%%%%%%%%%%%%%%%%%%%%% 
% in einem figure environment mit caption
   \calc_figscale{20}
    \begin{figure}[htb]
    \centerline{\input{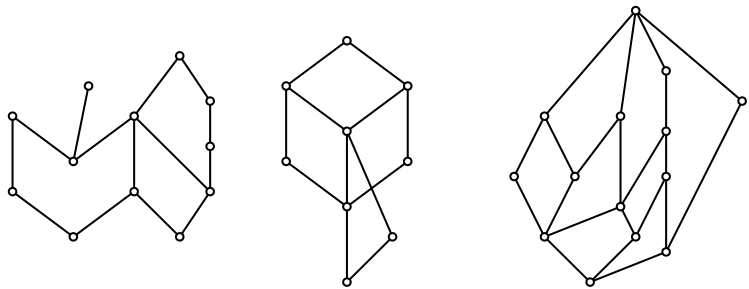}}
    \caption{Diagrams of a planar poset of dimension 3 (left), a
  non-planar lattice (middle), and a planar lattice (right).\label{fig:planar-posets}}
    \end{figure}
    
%%%%%%%%%%%%%%%%%%%%%%%%%%%%%%%%%%%%%%%%%%%%%%%%%%%%%%%%%%%%%%%%%%%%%%

For the reverse direction in the theorem, let $P$ be a lattice of
dimension at most~$2$.  Let $L_1$ and $L_2$ be linear orders on $X$ so
that $P=L_1\cap L_2$.  For each $x\in X$, and each $i=1,2$, let $x_i$
denote the height of $x$ in $L_i$.  Then a planar diagram of $P$ is
obtained by locating each $x\in X$ at the point in the plane with
integer coordinates $(x_1,x_2)$ and joining points $x$ and $y$ with a
straight line segment when one of $x$ and $y$ covers the other in
$P$. A pair of crossing edges in this drawing would violate the
lattice property, indeed if $x\prec y$ and $x'\prec y'$ are two covers
whose edges cross, then $x\leq y'$ and $x'\leq y$ whence $x$ and $x'$
have no unique least upper bound.
\smallskip

%%%%%%%%%%%%%%%%%%%%%%%%%%%%%%%%%%%%%%%%%%%%%%%%%%%%%%%%%%%%%%%%%%%%%% 
% in einem figure environment mit caption
   \calc_figscale{20}
    \begin{figure}[htb]
    \centerline{\input{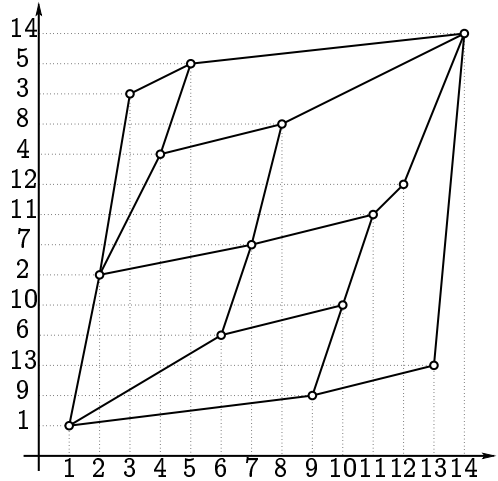}}
    \caption{The planar lattice from
  \cref{fig:planar-posets} with a realizer $L_1,L_2$.\label{fig:planar-lattice}}
    \end{figure}
    
%%%%%%%%%%%%%%%%%%%%%%%%%%%%%%%%%%%%%%%%%%%%%%%%%%%%%%%%%%%%%%%%%%%%%%
A planar digraph $D$ with a unique sink and source, both of them 
on the outer face, and no transitive edges is the digraph
of a planar lattice. Hence, the above discussion
directly implies the following classical result.

\begin{prop}
  A planar digraph $D$ on $n$ vertices with a unique sink and source
  on the outer face and no transitive edges has an upward
  drawing on an $(n-1)\times(n-1)$ grid. 
\end{prop}
To the best of our knowledge the area problem for diagrams of general
planar posets is open.

In this paper we will, henceforth, use the terms 2-dimensional poset
and planar lattice respectively to refer to a poset $P=(X,<)$ together
with a fixed ordered realizer $[L_1,L_2]$.  In the case of the
lattice, fixing the realizer can be interpreted as fixing a plane
drawing of the diagram. By fixing the realizer of $P$ we also have a
well-defined \term{primary conjugate}, this is the poset $Q$ on $X$
with realizer $[L_1,\ovl{L_2}]$, where $\ovl{L_2}$ is the reverse of
$L_2$.  Define the \term{left of relation} on $X$ such that $x$ is
left of $y$ if and only of $x=y$ or $x$ and $y$ are incomparable in
$P$ and $x<y$ in $Q$.

%%%%%%%%%%%%%%%%%%%%%%%%%%%%%%%%%%%%%%%%%%%%%%%%%%%%%%
\subsection{Orthogonal partitions of posets}
\label{ssec:GKtheory}

Let $P$ be a finite poset, Dilworth's theorem states that \textit{the
  maximum size of an antichain equals the minimum number of chains
  partitioning the elements of $P$}.

Greene and Kleitman \cite{GK} found a nice generalization of
Dilworth's result. Define a \term{$k$-antichain} to be a family of
$k$ pairwise disjoint antichains.

\begin{thm}\label{theo1}
For any partially ordered set $P$ and any positive integer $k$
$$\max\sum_{A\in \A}|A| = \min\sum_{C \in \C} \min(|C|,k) $$
where the maximum is taken over all $k$-antichains $\A$ and the minimum
over all chain partitions $\C$ of $P$.
\end{thm}

Greene \cite{G} stated the dual of this theorem. Let a
\term{$\ell$-chain} be a family of $\ell$ pairwise disjoint
chains.

\begin{thm}\label{theo3}
For any partially ordered set $P$ and any positive integer $\ell$
$$\max \sum_{C\in	\C}|C| = \min\sum_{A \in \A} \min(|A|,\ell)$$
where the maximum is taken over all $\ell$-chains $\C$ and the minimum
over all antichain partitions $\A$ of $P$.
\end{thm}

A further theorem of Greene \cite{G} can be interpreted as a generalization of
the Robinson-Schensted correspondence and its interpretation given by Greene
\cite{G/ST}.

To a partially ordered set $P$ with $n$ elements there is an associated 
partition $\lambda$ of $n$, such that for the Ferrer's diagram $G(P)$
corresponding to $\lambda$ we get:

\begin{thm}\label{theo5}
  The number of squares in the $\ell$ longest columns of $G(P)$ equals
  the maximal number of elements covered by an $\ell$-chain of $P$ and
  the number of squares in the $k$ longest rows of $G(P)$ equals the
  maximal number of elements covered by a $k$-antichain.
\end{thm}

\begin{figure}[htb]
    \centerline{\input{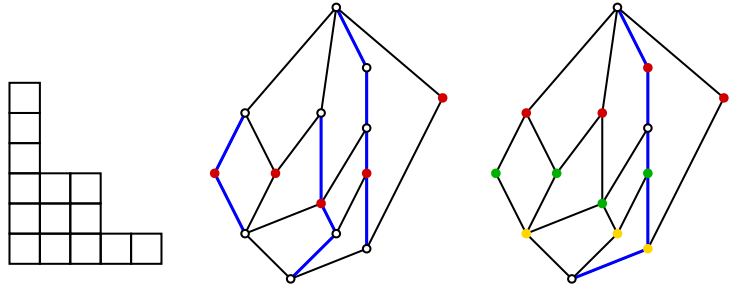}}
    \caption{The Ferrer's shape of the lattice $L$ from
  \cref{fig:planar-lattice} together with two orthogonal pairs of $L$
  corresponding to the boundary points $(1,3)$ and $(3,1)$ of $G(L)$;
  chains of $\C$ are blue, antichains of $\A$ are red, green, and yellow.\label{fig:gk-example}}
    \end{figure}

Figure~\ref{fig:gk-example} shows an example, in this case the
Ferrer's diagram $G(P)$ corresponds to the partition
$6+3+3+1+1\models 14$.  Several proofs of Greene's results are known,
e.g.~\cite{FOM},\cite{F}, and \cite{S}. For a not so recent, but at
its time comprehensive survey we recommend \cite{W}.

The approach taken by Andr\'{a}s Frank \cite{F} is particularly elegant.
Following Frank we call a chain family $\C$ and an antichain family $\A$ an 
\term{orthogonal pair} iff
\begin{enumerate}
\item \qquad$\displaystyle P = \Bigl(\bigcup_{A\in\A} A\Bigr)
\ \cup\ \Bigl(\bigcup_{C\in\C} C\Bigr) $, \ \textrm{and}
\item \qquad$|A\cap C| = 1 $\quad for all $A \in \A,\ C\in \C$.
\end{enumerate}
If $\C$ is orthogonal to a $k$-antichain $\A$ and $\C^{+}$ is obtained from
$\C$ by adding the rest of $P$ as singletons, then
$$
\sum_{A\in\A} |A| = \sum_{C\in\C^{+}} \sum_{A\in\A} 
|A \cap C| = \sum_{C\in\C^{+}} \min(|C|,k).
$$
Thus $\C^{+}$ is a $k$ optimal chain partition in the sense of Theorem~\ref{theo1}. Similarly an $\ell$
optimal antichain partition in the sense of Theorem~\ref{theo3} can be obtained from an orthogonal pair
$\A,\C$ where $\C$ is an $\ell$-chain.

Using the minimum cost flow algorithm of Ford and Fulkerson \cite{FF},
Frank proved the existence of a sequence of orthogonal chain and
antichain families. This sequence is rich enough to allow the
derivation of the whole theory. The sequence consists of an orthogonal
pair for every point from the boundary of $G(P)$.  With the point
$(k,\ell)$ from the boundary of $G(P)$ we get an orthogonal pair
$\A,\C$ such that $\A$ is a $k$-antichain and $\C$ an $\ell$-chain,
see Figure~\ref{fig:gk-example}. Since $G(P)$ is the Ferrer's diagram
of a partition of $n$ we can find a point $(k,\ell)$ on the boundary
of $G(P)$ with $k+\ell \leq \sqrt{2n}-1$ (This is because every
Ferrer's shape of a partition of $m$ which contains no point $(x,y)$
with $x+y \leq s$ on the boundary contains the shape of the partition
$(1,2,\ldots,s+1)$.  From $m \geq {s+2 \choose 2}$ we get
$s+1 < \sqrt{2m})$.

%%%%%%%%%%%%%%%%%%%%%%%%%%%%%%%%%%%%%%%%%%%%%%%%%%%%%%%%%%%%%%%%%%%%%% 
% in einem figure environment mit caption
   \calc_figscale{20}
    
%%%%%%%%%%%%%%%%%%%%%%%%%%%%%%%%%%%%%%%%%%%%%%%%%%%%%%%%%%%%%%%%%%%%%%

We will use the following corollary of the theory:
\begin{cor}\label{cor:ortho-pair}
Let $P=(X,<)$ be a partial order on $n$ elements, then there is an
orthogonal pair $\A,\C$ where $\A$ is a $k$-antichain and
$\C$ an $\ell$-chain and $k+\ell \leq \sqrt{2n}-1$.
\end{cor}

For our application we will need some additional structure on the
antichains and chains of an orthogonal pair $\A,\C$.

The \term{canonical
antichain partition} of a poset $P=(X,<)$ is constructed by recursively
removing all minimal elements from $P$ and make them one of the
antichains of the partition. More explicitely
$A_1=\Min(X)$ and  $A_j = \Min\big(X\setminus \bigcup\{A_i : 1\leq i < j\}\big)$
for $j>1$. Note that by definition
for each element $y\in A_j$ with $j>1$ there is some $x\in A_{j-1}$
with $x<y$. Due to this property there is a chain of $h$ elements in
$P$ if the canonical antichain partition consists of $h$ non-empty
antichains. This in essence is the dual of Dilworth's theorem, i.e.,
the statement: \textit{the
  maximal size of a chain equals the minimal number of antichains
  partitioning the elements of $P$}.

\begin{lemma}\label{lem:cap-in-op}
  Let $\A,\C$ be an orthogonal pair of $P=(X,<)$ and let $P_\A$ be the
  order induced by $P$ on the set $X_\A = \bigcup\{A: A \in \A \}$.
  If $\A'$ is the canonical antichain partition of $P_\A$, then
  $\A',\C$ is again an orthogonal pair of $P$
\end{lemma}

\Proof Let $\A$ be the family $A_1,\ldots,A_k$. Starting with this
family we will change the antichains in the family
while maintaining the invariant that
the family of antichains together with $\C$ forms an
orthogonal pair. At the end of the process the family of antichains
will be the canonical antichain partition of $P_\A$.

The first phase of changes is the uncrossing phase. We iteratively
choose two antichains $A_i,A_j$ with $i<j$ from the present
family and let
$B_i = \{y\in A_i : \textrm{ there is an } x\in A_j \textrm{ with }
x<y\}$ and
$B_j = \{x\in A_j : \textrm{ there is a } y\in A_i \textrm{ with }
x<y\}$. Define $A'_i = A_i-B_i+B_j$ and $A'_j = A_j-B_j+B_i$. It
is easy to see that $A'_i$ and $A'_j$ are antichains and that the
family obtained by replacing $A_i,A_j$ by $A'_i,A'_j$ is
orthogonal to $\C$. This results in a family of $k$ antichains such
that if $i<j$ and $x\in A_i$ and $y\in A_j$ are comparable, then $x<y$.

The second phase is the push-down phase.  We iteratively choose
$i\in[k-1]$ and let
$B = \{y\in A_{i+1} : \textrm{ there is no } x\in A_i \textrm{ with }
x<y\}$ and define $A'_{i+1} = A_{i+1}-B$ and $A'_i = A_i+B$.  It is
again easy to see that $A'_i$ and $A'_{i+1}$ are antichains and that
the family obtained by replacing $A_i,A_{i+1}$ by $A'_i,A'_{i+1}$ is
orthogonal to $\C$. This results in a family of $k$ antichains such
that if $y\in A_{i+1}$, then there is an $x\in A_i$ with
$x<y$. This implies that $A_j = \Min(X_\A\setminus \bigcup\{A_i :
1\leq i < j\})$, whence the family is the canonical antichain partition.
\qed
%% ------------------------------------------

Let $P=(X,<)$ be a 2-dimensional poset with realizer $[L_1,{L_2}]$ and recall
that the primary conjugate has realizer $[L_1,\ovl{L_2}]$.  The order $Q$
corresponds to a transitive relation on the complement of the comparability
graph of $P$, in particular chains of $P$ and antichains of $Q$ are in
bijection.

The canonical antichain partition of $Q$ yields the 
\term{canonical chain partition} of $P$. The canonical chain
partition $C_1,C_2,\ldots,C_w$ of $P$ can be characterized by the
property that  for each $1\leq i < j\leq w$ and each element $y\in C_j$
there is some $x\in C_{i}$ with $x\; ||\; y$ and in $L_1$
element $x$ comes before $y$. In particular $C_1$ is a maximal chain of $P$.

Let $\A,\C$ be an orthogonal pair of the 2-dimensional $P=(X,<)$.
Applying the proof of \cref{lem:cap-in-op} to the orthogonal
pair $\C,\A$ of $Q$ we obtain:

\begin{lemma}\label{lem:ccp-in-op}
  Let $\A,\C$ be an orthogonal pair of $P=(X,<)$ and let $P_\C$ be the
  order induced by $P$ on the set
  $X_\C = \bigcup\{C : C\in \C\}$.  If $\C'$ is the canonical
  chain partition of $P_\C$, then $\C',\A$ is again an orthogonal
  pair of $P$
\end{lemma}

In a context where edges of the diagram are of interest, it is convenient to
work with maximal chains. The canonical chain partition $C_1,C_2,\ldots,C_w$
of a 2-dimensional $P$ induces a \term{canonical chain cover} of $P$ which
consists of maximal chains.  With chain $C_i$ associate a chain $\mC_i$ which
is obtained by successively adding to $C_i$ all compatible elements of
$C_{i-1},C_{i-2},\ldots$ in this order. Alternatively $\mC_i$ can be described
by looking at the conjugate of $P$ with realizer $[\ovl{L_1},{L_2}]$ (this is
the dual of the primary conjugate $Q$), and defining $\mC_i$ as the first chain
in the canonical chain partition of the order induced by
$\bigcup\{C_j : 1\leq j \leq i\})$. The maximality of $\mC_i$ follows from the
characterization of the canonical chain partition given above.

%%%%%%%%%%%%%%%%%%%%%%%%%%%%%%%%%%%%%%%%%%%%%%%%%%%%%%
\section{Drawing Planar Lattices on Few Lines}
\label{sec:plaLat}
%%%%%%%%%%%%%%%%%%%%%%%%%%%%%%%%%%%%%%%%%%%%%%%%%%%%%%

In this section we prove that planar lattices with $n$ elements have
a straight-line diagram with all vertices on a set of $\sqrt{2n}-1$
horizontal and vertical lines. The following proposition covers the
case where the lattice has an antichain partition of small size.
We assume that a planar lattice is given with a realizer $[L_1,L_2]$
and, hence, with a fixed plane drawing of its diagram.

\begin{prop}\label{prop:prescribed-h}
  For any planar lattice $L=(X,<)$ with an extension $h:X \to \RR$ of
  $L$ there is a plane straight-line drawing $\Gamma$ of the diagram
  $D_L$ of $L$ such that each element $x\in X$ is represented by a
  point with $y$-coordinate $h(x)$. Additionally all elements of the
  left boundary chain of $D_L$ are aligned vertically in the drawing.
\end{prop}

\Proof
Let $C_1,C_2,\ldots,C_w$ be the canonical chain partition and
$\mC_1,\mC_2,\ldots,\mC_w$ be the corresponding canonical chain cover.
Define $S_i$ as the suborder of $L$ induced by
$\bigcup\{C_j : 1\leq j \leq i \}$ and note that $S_i$ is a sublattice
of $L$ with left boundary chain $C_1=\mC_1$ and right boundary chain
$\mC_i$.

Embed the elements of $C_1$ on a vertical line $\bg_1$ (e.g.~the line
$y=0$) with points as prescribed by $h$. This is a drawing $\Gamma_1$
of $S_1$. Suppose that a drawing $\Gamma_i$ of the diagram $S_i$ is
constructed. The right boundary path $\gamma_i$ of $\Gamma_i$ is a polygonal
$y$-monotone path. Embed the elements of $C_{i+1}$ on a vertical line
$\bg_{i+1}$ with points as prescribed by $h$. We need a position for
$\bg_{i+1}$ to the right of $\gamma_i$ such that all the diagram edges
connecting $C_{i+1}$ to $\mC_i$ can be inserted to obtain a crossing
free drawing $\Gamma_{i+1}$ of the diagram of $S_{i+1}$.

Let $E_i$ be the set of diagram edges connecting $C_{i+1}$ to
$\mC_i$. For each $e\in E_i$ there are points $p\in \gamma_i$ and
$q\in \bg_{i+1}$ representing the endpoints.  Let $K_p$ be an open
cone with apex $p$ which intersects $\gamma_i$ only at $p$ and
contains a horizontal ray to the right.  Let $b_e$ be the minimal
horizontal distance of $\gamma_i$ and $\bg_{i+1}$ such that
$q\in K_p$. Let $\beta = \max( b_e : e\in E_i )$.  If we place
$\gamma_i$ and $\bg_{i+1}$ at horizontal distance $\beta$, then the
edges of $E_i$ can be drawn such that they do not interfere (introduce
crossings) with $\gamma_i$. We claim that there is no crossing of
edges of $E_i$. Let $(p,q)$ and $(p',q')$ be two drawn edges from
$E_i$. Since they are edges of a planar diagram and have endpoints on
two chains, we know, that $h(p)\leq h(p')$ implies $h(q)\leq h(q')$.
Edge $(p,q)$ is drawn in the cone $K_p$. If $(p',q')$ intersects the
edge and $p'$ is above $p$ on $\gamma_i$, then $q'$ has to be below
$q$ on $\bg_{i+1}$, in contradiction to $h(q)\leq h(q')$. Hence we
have a planar drawing $\Gamma_{i+1}$ of the diagram of $S_{i+1}$.
With induction we obtain the drawing $\Gamma=\Gamma_w$ of $D_L$.  \qed

\begin{thm}\label{thm:lat-few}
  For every planar lattice $L=(X,<)$ with $|X|=n$, there is a plane
  straight-line drawing of the diagram such that the elements are
  represented by points on a set of at most $\sqrt{2n}-1$
  lines. Additionally
\Bitem each of the lines is either horizontal or
  vertical,
\Bitem each crossing point of a horizontal and a vertical
  line hosts an element of~$X$.
\end{thm}

\Proof 
Let $\A,\C$ be an orthogonal pair of $L$ such that
$\A$ is a $k$-antichain, $\C$ an $\ell$-chain, and $k+\ell \leq \sqrt{2n}-1$.
It follows from \cref{cor:ortho-pair} that such a pair exists.

Since $L$ has a fixed ordered realizer $[L_1,L_2]$, we can apply
\cref{lem:cap-in-op} to $\A$ and \cref{lem:ccp-in-op} to $\C$ to get an
orthogonal pair $(A_1,\ldots,A_k),(C_1,\ldots,C_\ell)$ where the antichain
family and the chain family are both canonical.  Fix an extension $h:X\to \RR$
of $L$ with the property that $h(x) = i$ for all $x\in A_i$.

In the following we will construct a drawing
$\Gamma$ of $D_L$ such that each element $x\in X$ is represented by a
point with $y$-coordinate $h(x)$, and in addition all elements of
chain $C_i$ lie on a common vertical line $\bg_i$ for $1\leq i\leq \ell$ .
By Property~1 of orthogonal pairs, for each $x\in X$ there is an
$i$ such that $x\in A_i$ or a $j$ such that $x\in C_j$ or both.
Therfore, $\Gamma$ will be a drawing such that the $k$ horizontal lines
$y=i$ with $i=1,\ldots,k$ together with the $\ell$ vertical
lines $\bg_j$ with $j=1,\ldots,\ell$ cover all the elements of $X$.
Property 2 of orthogonal pairs
implies the second extra property mentioned in the theorem.

If the number $\ell$ of chains is zero, then we get a drawing
$\Gamma$ with all the necessary properties from
\cref{prop:prescribed-h}. Now let $\ell > 0$.

The chain family $C_1,\ldots,C_\ell$ is the canonical chain partition
of the order induced on $X_\C = \bigcup\{C_i : i=1\ldots \ell\}$.
Let $\mC_1,\ldots,\mC_\ell$ be the corresponding canonical chain
covering of $X_\C$.

Let $X_i$ for $1\leq i\leq \ell$ be the set of all elements which are
left of some element of $\mC_i$ in $L$, and let $X_{\ell+1}=X$. Define
$S_i$ as the suborder of $L$ induced by $X_i$.  Also let
$Y_i = X_{i+1} - X_{i} + \mC_i$ and let $T_i$ be the suborder of $L$
induced by $Y_i$. Note the following properties of these sets and
orders:
\Bitem
$X_i\cap C_j = \emptyset$ for $1\leq i<j\leq\ell$.
\Bitem
Each $S_i$ is a planar sublattice of $L$, its right boundary chain is
$\mC_i$.
\Bitem
$T_i$ is a planar sublattice of $L$.
\medskip

\ni
A drawing $\Gamma_1$ of $S_1$ with the right boundary chain being
aligned vertically is obtained by applying \cref{prop:prescribed-h}
to the vertical reflection of the diagram $D_L[X_1]$ and
reflecting the resulting drawing vertically.

We construct the drawing $\Gamma$ of $D_L$ in phases.  
In phase $i$ we aim for a drawing $\Gamma_{i+1}$ of $S_{i+1}$ extending
the given drawing $\Gamma_{i}$ of $S_i$, i.e., we need to construct a
drawing $\Lambda_i$ of $T_i$ such that
\Item{(1)} The left boundary chain of
$\Lambda_i$ matches the right boundary chain of $\Gamma_{i}$.
\Item{(2)} In $\Lambda_i$ all elements of $C_{i+1}$ are on a common
vertical line $\bg_{i+1}$.
\smallskip

\ni
The construction of $\Lambda_i$ is done in three stages. First we
extend $\mC_i$ to the right by adding `ears'. Then we extend $\mC_{i+1}$
to the left by adding `ears'. Finally we show that the left and the right
part obtained from the first two stages can be combined to yield the
drawing $\Lambda_i$.

To avoid extensive use of indices let $Y=Y_i$, $T=T_i$, $\mC=\mC_i$,
and let $\gamma$ be a copy of the $y$-monotone polygonal
right boundary of $\Gamma_i$, i.e., $\gamma$ is a drawing of~$C$.
We initialize $\Lambda'=\gamma$.

A \term{left ear} of $T$ is a face $F$ in the diagram $D_L[Y]$ of $T$
such that the left boundary of $F$ is a subchain of the left boundary
chain $\mC$ of $D_L[Y]$. The ear is \textit{feasible} if the right
boundary chain contains no element of $C_{i+1}$. Given a feasible ear
we use the method from the proof of \cref{prop:prescribed-h} to add
$F$ to $\gamma$.  We represent the right boundary
$z_0 < z_1 < \ldots < z_{l}$ excluding $z_0$ and $z_l$ of $F$ on a
vertical line $\bg$ by points $q_1,\ldots,q_{l-1}$ with
$y$-coordinates as prescribed by $h$.  The points $q_0$ and $q_{l}$
representing $z_0$ and $z_l$ respectively are already represented on
$\gamma$. Then we place $\bg$ at some distance $\beta$ to the right of
$\gamma$. The value of $\beta$ has to be chosen large enough to ensure
that edges $q_0,q_1$ and $q_{l-1},q_l$ are drawn such that they do not
interfere with $\gamma$.  Let $\Lambda'$ be the drawing augmented by
the polygonal path $q_0,q_1,\ldots,q_{l-1},q_l$ and let $\mC$ again
refer to the right boundary chain $\gamma$ of $\Lambda'$. Delete all
elements of the left boundary of $F$ except $z_0$ and $z_l$ from $Y$
and $T$. This \textit{shelling of a left ear} from $T$ is iterated
until there remains no left feasible ear. Upon stopping we have a
drawing $\Lambda'$ which can be glued to the right side of
$\Gamma_i$.  Let $\gamma'$ be the right boundary chain of $\Lambda'$.

Now let $C=C_{i+1}$.  Initialize a new drawing $\Lambda''$ by placing the
elements of $C$ on a vertical line $\bg$ at the heights prescribed by $h$ and
connect consecutive ones by an edge whenever the order relation is indeed a
cover relation of $L$. The initial drawing may thus be disconnected and if so
this will remain the case throughout this stage. We now consider right ears
from $T$. A right ear of $T$ corresponding to a face $F$ is feasible
if the left boundary chain of $F$ contains no element of $\gamma'$. The left
boundary chain of a feasible ear can be drawn as a $y$-monotone polygonal
chain left of the left boundary $\gamma''$ of $\Lambda''$. Update $\gamma''$
to be the new left boundary of the augmented $\Lambda''$ and remove the
elements of the ear from $Y$ and $T$.  The shelling of right ears from $T$ is
iterated until there remains no feasible right ear. Note that $\gamma''$ is
$y$-monotone but it may consist of several components.

In the final stage we have to combine the drawings
$\Lambda'$, $\Lambda''$ into a single drawing. This is done by drawing the edges and
chains which remain in $T$ between the two boundary chains as
straight segments between $\gamma'$ and $\gamma''$.  This will be
possible because we can shift $\gamma'$ and $\gamma''$ as far apart
horizontally as necessary.

First we draw all the edges connecting the two chains.
Let $E$ be the set of edges connecting the left and right boundary
chains of $T$. For each $e\in E$ there are points $p\in \gamma'$ and
$q\in \gamma''$ representing the endpoints.  Let $K_p$ be an open cone with
apex $p$ which intersects $\gamma'$ only at $p$ and contains a
horizontal ray to the right and let $K_q$ be an open cone with apex $q$
which intersects $\gamma''$ only at $q$ and contains a horizontal ray
to the left. Let $b_e$ be the minimal horizontal distance of $\gamma'$
and $\gamma''$ such that $p\in K_q$ and $q\in K_p$. Let
$\beta = \max( b_e : e\in E )$.  If we place $\gamma'$ and $\gamma''$
at horizontal distance $\beta$, then the edges of $E$ can be drawn
such that they do not interfere (introduce crossings) with $\gamma'$
and $\gamma''$. We claim that there is no crossing of edges of
$E$. Let $(p,q)$ and $(p',q')$ be two drawn edges from $E$. Since they
are diagram edges with endpoints on two chains we know that
$h(p)\leq h(p')$ implies $h(q)\leq h(q')$.  Edge $(p,q)$ is drawn in
$K_p\cap K_q$. If $(p',q')$ intersects the edge and $p'$ is above $p$
on $\gamma'$, then $q'$ has to be below $q$ on $\gamma''$, in
contradiction to $h(q)\leq h(q')$. Placing $\Lambda'$ and $\Lambda''$
such that $\beta$ is the distance between their outer chains and
drawing the edges of $E$ yields a drawing $\Lambda$ of some lattice.
An important feature of the drawing is that if
we move the two subdrawings $\Lambda'$ and $\Lambda''$ further
apart the drawing keeps the needed properties, i.e., the height of
elements remains unaltered, vertices of a chain which should be
vertically aligned remain vertically aligned, and the drawing is
crossing-free.

Now assume that $T$ contains elements which are not represented in $\Lambda$.
Let $B$ be a connected component of such elements where connectivity is with
respect to $D_L$. All the elements of $B$ have to be placed in a face $F_B$ of
$\Lambda$. Let $\delta'$ and $\delta''$ be the left and right boundary of
$F_B$.

In the following we will repeat the choice of a component $B$ and a
chain $C$ from $B$ which is to be drawn in the corresponding face
$F_B$ of $\Lambda$ such that the minimum and the maximum of $C$ have
connecting edges to the two sides of the boundary of $F_B$. Let us
consider the case that in $D_L$ the maximum of $C$ has an outgoing edge to
an element which is represented by a point $p\in \delta'$ and the
minimum of $C$ has an incoming edge from an element
represented by $q\in\delta''$. We represent the elements of $C$ as
points on the prescribed heights on a line segment $\zeta$ with
endpoints $p$ and~$q$. It may become necessary to stretch the face
horizontally to be able to place $C$. In this case we stretch the
whole drawing between $\gamma'$ and $\gamma''$ with a uniform stretch
factor. There may be additional edges between elements of $C$ and
elements on $\delta'$ and $\delta''$. They can also be drawn without
crossing when the distance of $\delta'$ and $\delta''$ exceeds some
value $b$.

Stretching the
whole drawing between $\gamma'$ and $\gamma''$ allows us to draw
the segment $\zeta$ and additional edges inside of $F_B$ because
of the following invariant.

\Bitem For each face $F$ of the drawing $\Lambda$ and two points $x$ and $y$
from the boundary of $F$ it holds that: if the segment $x,y$
is not in the interior of $F$, then the parts of the boundary obstructing the
segment $x,y$ belong to $\gamma'$ or $\gamma''$.
\smallskip

\ni When including a chain $C$ in the drawing $\Lambda$, we place the elements
of $C$ at the prescribed heights on a common line segment $\zeta$. This
ensures that each new element contributes convex corners in all incident
faces. Hence, new elements can not obstruct a visibility within a
face. Therefore, obstructing corners correspond to elements of $\gamma'$ or
$\gamma''$ and the invariant holds.

Now consider the case where maximum and minimum of the chain $C$
connect to two elements $p$ and $q$ on the same side of $F$.  Since
$\gamma'$ and $\gamma''$ do not admit ear extensions we know that not
both of $p$ and $q$ belong to one of $\gamma'$ and $\gamma''$.  If the
segment from $p$ to $q$ is obstructed, then the invariant ensures that
with sufficient horizontal stretch the segment $\zeta$ connecting $p$
and $q$ will be inside $F$. Hence, chain $C$ can be drawn and
$\Lambda$ can be extended.

When there remains no component $B$ containing a chain $C$ which can be
included in the drawing using the above strategy, then either all
elements of $Y$ are drawn or we have the
following: every component $B$ only connects to elements of a
line segment $\zeta_B$.

In this situation $B$ is kind of a big ear over $\zeta_B$.
We next describe how to draw $B$, but note, that doing this we
will not maintain or need the invariant. 

By construction all elements of $\zeta_B$
belong to a common chain $C_B$. Consider the union
$B + C_B$ and note that this is a planar lattice $L_B$,
moreover, $C_B$ is either the left or the right boundary chain
of $L_B$. Assume that $C_B$ is the left boundary
chain of $L_B$. Now use \cref{prop:prescribed-h} to get a drawing
$\Lambda_B$ of $L_B$ with $C_B$ aligned vertically. Using an
affine transformation we can map $\Lambda_B$ into $\Lambda$
such that the line containing $C_B$ in $\Lambda_B$ is mapped
to the line supporting the segment $\zeta_B$. Since elements of $C_B$
are at their prescribed heights their representing points in
$\Lambda_B$ are mapped to the representing points of $\Lambda$. 
The affine map also has to compress $\Lambda_B$ horizontally so that it
is placed in a narrow strip on the right side of $\zeta_B$. This
strip can be chosen narrow enough to make sure that 
all of $B$ is mapped to the face of $\Lambda$
where it belongs.

Glueing the drawings $\Lambda'$ with $\Lambda$
at the polygonal path $\gamma'$ and $\Lambda$ with $\Lambda''$
at $\gamma''$ (a $y$-monotone collection of paths) yields
a drawing $\Lambda_i$ of $T_i$. The drawing $\Lambda_i$ can be
glued to $\Gamma_i$ to form a drawing $\Gamma_{i+1}$ of $S_{i+1}$.
Eventually the drawing $\Gamma_\ell$ will be constructed.
From there the drawing $\Gamma=\Gamma_{\ell+1}$ is obtained by adding
some left ears.
\qed

%%%%%%%%%%%%%%%%%%%%%%%%%%%%%%%%%%%%%%%%%%%%%%%%%%%%%%
\section{Transversal Structures on Few Lines}
\label{sec:4ct}

\begin{thm}\label{thm:draw-ts}
For every internally 4-connected inner
triangulation of a 4-gon $G=(V,E)$ with $n$ vertices
there is a planar straight line drawing  such that the vertices are
  represented by points on a set of at most $\sqrt{2n}-1$
  lines. Additionally
\Bitem each of the lines is either horizontal or
  vertical,
\Bitem each crossing point of a horizontal and a vertical
  line hosts a vertex.
\end{thm}

\Proof
Fix a transversal structure of $G$ and consider the red graph
$G_R=(V,E_R^+)$. From \cref{prop:bipolar} and (R) we know that
$G_R$ is bipolar and transitively reduced. This implies that
there is a planar lattice $L=(V,<)$ such that a diagram of $L$ is an
upward drawing of $G_R$. The relation $<$ is defined as $v < v'$ 
if and only if there is a directed path from $v$ to $v'$ in $G_R$.

We would like to use \cref{thm:lat-few} to draw $G_R$ on $\sqrt{2n}-1$
lines and then include the blue edges of the transversal structure in
the drawing. This, however, may yield crossings.
Instead we go through the proof of \cref{thm:lat-few} and
include blue edges while constructing the drawing of the red graph.

When adding a left feasible ear, i.e., when adding the right boundary
of a face $F$, we draw all the blue edges
corresponding to the face $F$. If $e$ has to connect $p\in \gamma$
and $q \in \bg$ define $b_e$ as the minimal
horizontal distance of $\gamma$ and $\bg$ such that
$q\in K_p$. When placing $\bg$ at a distance $\beta$ from $\gamma$
which exceeds all the values $b_e$, the blue edges can be drawn
crossing free.
When adding a right feasible ear the situation is symmetric.

Now let us consider the stage where
a left and right drawing $\Lambda'$ and $\Lambda''$
with boundary chains $\gamma'$ and $\gamma''$ have to be combined.
When drawing edges connecting the two chains we include the set of all
blue edges with one end on $\gamma'$ and one on $\gamma''$.
Then we complete the combination on the basis of the red edges.
Only in the `bad' case we have to be careful.
First, when drawing $L_B$ using \cref{prop:prescribed-h} we
also include the blue edges in the drawing. This only requires to
choose the distances $\beta$ as maxima over larger sets of values
$b_e$. Second, when placing the drawing $\Lambda_B$ in a narrow strip on the
side of $\zeta_B$ we have to be carefull that this does not obstruct a
visibility from the left side of the face to the right side.
Finally, all the remaining blue edges have to be drawn in the faces
between $\gamma'$ and $\gamma''$. Due to the invariant this is
possible if we stretch the drawing between the two chains
sufficiently.
\qed

%%%%%%%%%%%%%%%%%%%%%%%%%%%%%%%%%%%%%%%%%%%%%%%%%%%%%%%%%%%%%%%%%%%%%% 
% in einem figure environment mit caption
   \calc_figscale{30}
    \begin{figure}[htb]
    \centerline{\input{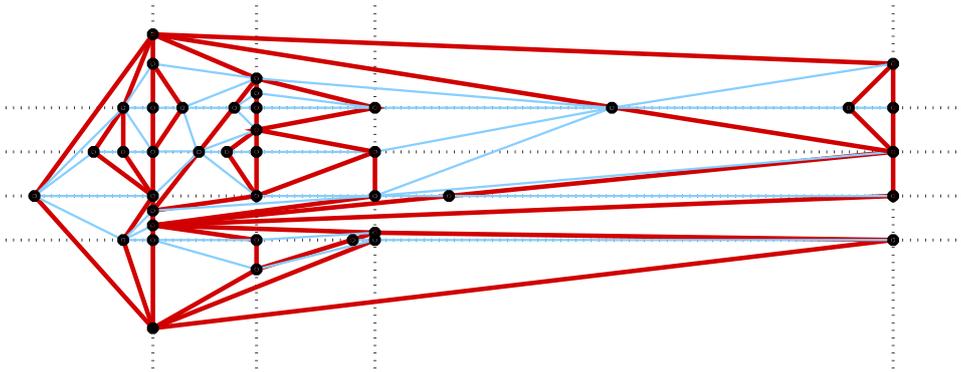}}
    \caption{A partially drawn transversal structure.
  The figure shows a drawing of $\Gamma_4$, these are the vertices
  left of some element in $\mC_4$ together with the induced edges.\label{fig:big-example}}
    \end{figure}
    
%%%%%%%%%%%%%%%%%%%%%%%%%%%%%%%%%%%%%%%%%%%%%%%%%%%%%%%%%%%%%%%%%%%%%%

It remains to see how \cref{thm:main} follows from \cref{thm:draw-ts}.
Let $G$ be a 4-connected triangulation and let $G'$ be obtained from
$G$ by deleting one of the outer edges. Now $G'$ is an internally
4-connected inner triangulation of a 4-gon.  Label the outer vertices
of $G'$ such that the deleted edge is the edge $s,t$.  Slightly
stretching \cref{thm:draw-ts} we prescribe $h(s)=-\infty$ and
$h(t)=\infty$, this yields a planar straight-line drawing $\Gamma$ of
$G'$ such that the vertices except $s$ and~$t$ are represented by
points on a set of at most $\sqrt{2n}-1$ lines and the edges
connecting to $s$ and~$t$ are vertical rays.  Moreover with every edge
$v,s$ or $v,t$ there is an open cone $K$ containing the vertical ray,
such that the point representing $v$ is the apex of $K$ and this is
the only vertex contained in $K$. Now let $\bg$ be a vertical line
which is disjoint from $\Gamma$. On $\bg$ we find a point $p_s$ which
is contained in all the upward cones and a point $p_t$ contained in
all the downward cones.  Taking $p_s$ and $p_t$ as representatives for
$s$ and $t$ we can tilt the rays and make them finite edges ending in
$p_s$ and $p_t$ respectively, and in addition draw the edge $p_s,p_t$.
\medskip

We conclude with a remark and two open problems.
\Bitem Our results are constructive and can be complemented with
algorithms running in polynomial time.
\Bitem Is $\pi(G) \in O(\sqrt{n})$ for every planar graph $G$ on $n$ vertices?
\Bitem What size of a grid is needed for drawings of 4-connected plane
graphs on $O(\sqrt{n})$ lines?
\medskip

%%%%%%%%%%%%%%%%%%%%%%%%%%%%%%%%%%%%%%%%%%%%%%%%%%%%%%
\subsection*{Acknowledgments}
Work on this problem began at the 2018 Bertinoro Workshop of Graph
Drawing. I thank the organizers of the event for making this
possible. Special thanks go to Pavel Valtr, Alex Pilz and Torsten
Ueckerdt for helpful discussions.

% *************************************************************
% *************************************************************
\vskip-4mm\vbox{}
%%\small
\raggedright
% *************************************************************
\global\advance\baselineskip -.2pt
\advance\bibitemsep -2pt
% *************************************************************
\def\doiitem#1{\hbox{\footnotesize Doi:\;\url{#1}}}
\let\sc=\SC
\bibliography{lines-bib} 
\bibliographystyle{my-siam}

\end{document}